\documentclass[12pt]{article}

\usepackage{amsmath,amsfonts,amssymb,amsthm}
\usepackage{cite}
\usepackage{graphics}

\usepackage{epstopdf}

\usepackage{hyperref} 
\hypersetup{colorlinks=true, linkcolor=blue,  anchorcolor=blue,  
citecolor=blue, filecolor=blue, menucolor=blue, pagecolor=blue,  
urlcolor=blue,pdftitle={ACmain}}

\title{ \textbf{A connection between Tests for absolute convergence of infinite series, or how to be fair}}
\author{\emph{Victoria Rayskin}\\ \small{Optimum Solvers}\\ \small{vrayskin@gmail.com}}
\date{}

\theoremstyle{plain}
\newtheorem{theorem}{Theorem}[section]
\newtheorem{lemma}[theorem]{Lemma}
\newtheorem{corollary}[theorem]{Corollary}

\theoremstyle{definition}
\newtheorem{definition}[theorem]{Definition}
\newtheorem{example}[theorem]{Example}

\theoremstyle{remark}
\newtheorem{remark}{Remark}

\newcommand{\R}{\mathbb{R}}
\newcommand{\N}{\mathbb{N}}

\begin{document}

\maketitle

\begin{abstract}
The Ratio Test and the Root Test for absolute convergence/divergence of series of numbers $\sum_{n=0}^{\infty}a_n$ are frequently discussed and proved independently in Calculus courses.

The Root Test is stronger (verifies convergence for more series) than the Ratio Test. This relation inspires introduction of some intermediate strength tests (stronger than the Ratio Test and weaker than the Root Test) that we call Power Mean Tests (they, in particular, include the Arithmetic Mean Test). We show the connection between the Root, the Power Mean and the Ratio Tests. We also note that all these tests are related to the test that we formulate and call Generalized $f$-mean Test (or Kolmogorov-Nagumo-de Finetti mean Test). We provide an example of an infinite series, where the Arithmetic Mean test is the test that should be used for convergence verification, because the Root and the Ratio Tests are not easy to apply.


We conclude this work with the statement emphasizing why it is ``fair'' to include the summarizing Corollary~\ref{connections} in our Calculus course.
\end{abstract}


\section{Introduction}\label{sec-intro}

In the study of infinite series convergence tests, we usually discuss the Root Test and the Ratio Test. Namely, let us consider the series of non-zero real numbers
\[
\sum_{n=0}^{\infty}a_n.
\]
Define 
\[
\rho_{\mbox{ratio}}=\lim_{n\to\infty}\frac{|a_{n+1}|}{|a_n|}
\]
and 
\[
\rho_{\mbox{root}}=\lim_{n\to\infty}\sqrt[n]{|a_n|}.
\]

The Ratio Test states that if $\rho_{\mbox{ratio}}<1$, the series is absolutely convergent. If $\rho_{\mbox{ratio}}>1$, the series is divergent. If $\rho_{\mbox{ratio}}=1$, the Ratio Test is inconclusive.

A similar result can be stated for the Root Test. If $\rho_{\mbox{root}}<1$, the series is absolutely convergent. If $\rho_{\mbox{root}}>1$, the series is divergent. If $\rho_{\mbox{root}}=1$, the Root Test is inconclusive.

The application of each of these tests is valuable, as some series are easier to test with the Root Test, and the Ratio Test is more suitable for the others. Root Test is stronger (Theorem~\ref{root-test}) than the Ratio Test (it verifies convergence for more series). We discuss here intermediate strength tests that we call Power Mean Tests for convergence and their further generalization, Generalized $f$-mean Test (Theorem \ref{generalized}):
\\
{\it
Let $\sum_{n=0}^{\infty} a_n$ be a series of non-zero real numbers, $I\subset \R$ be a finite or infinite interval, such that $\left|\frac{a_1}{a_0}\right|, \left|\frac{a_2}{a_1}\right|, ... \in I $.
Let $f$ be a continuous injective function, $f: I \to \R$, and assume that the limit
$$
\rho_{\mbox{generalized}}=\lim_{N\to\infty} f^{-1}\left(\frac{1}{N}\sum_{n=1}^{N}f\left|\frac{a_n}{a_{n-1}}\right|\right).
$$
exists.\\
If $\rho_{\mbox{ratio}}$ exists, then,  the value of the $\rho_{\mbox{generalized}}$ is independent of the choice of $f$. In particular, $\rho_{\mbox{generalized}}=\rho_{\mbox{ratio}}$.
}

When $f$ is a power functions, we obtain the Power Mean Test (Theorem~\ref{power}):
\\
{\it
Let $\sum_{n=0}^{\infty} a_n$ be a series of non-zero real numbers, let $p$ be a real positive number and suppose that 
\[
\rho_{\mbox{mean}} := \lim_{N\to\infty}\sqrt[p]{\frac{1}{N}\sum_{n=1}^{N}\left|\frac{a_n}{a_{n-1}}\right|^p}.
\]
If $\rho_{\mbox{mean}}<1$, then the series is convergent. If $\rho_{\mbox{ratio}}$ exists,  then $\rho_{\mbox{ratio}}=\rho_{\mbox{mean}}=\rho_{\mbox{root}}.$ 
}

If $p_1$ and $p_2$ are the two exponents in the Power Mean Test, and $p_1 < p_2$, then the Power Mean Test associated with $p_1$ is stronger than the Power Mean Test associated with $p_2$ (see Remark~\ref{p-comparison}).

When $p=1$ the Power Mean Test can be called Arithmetic Mean Test. 

We present an example of a series that is convergent according to the Arithmetic Mean Test. However, neither the Ratio Test, nor the Root Test is easy to use for this series. See Example~\ref{example}.

We also note that the Ratio, the Root and the Power Mean Tests are related to the Generalized $f$-mean Test (see Remark~\ref{instances}).

There are many examples (\cite{HWT}, \cite{C-U}, \cite{B}, \cite{IP}, \cite{T}, \cite{BB}) showing that the limit in the Root Test may exist while the limit in the Ratio Test does not. 
We provide simple proofs showing that the Root Test is stronger than the Mean Value Test and the latter is stronger than the Ratio Test. See Theorems~\ref{root-test} and \ref{power}. The benefits of these proofs are two-fold.

First, the proofs clarify that the radius of convergence of series of functions (and power series in particular) is well defined and can be calculated with either Root or Ratio Tests.

Second, Corollary~\ref{connections} illustrates that in the case of $\rho=1$ the inconclusiveness of one test implies inconclusiveness of some other tests.

\section{Convergence Tests}\label{proof}
First, we will recall the two well-known tests: the Ratio Test and the Root Test, and then we will state and prove the Theorem~\ref{root-test} that explains the relation between the two tests and shows that the root test is stronger. After that we introduce an intermediate strength test that we call Power Mean Test (Theorem~\ref{power}) and demonstrate its application with the Example~\ref{example}.
\begin{theorem}[The Ratio Test]
Let $\sum a_n$ be a series of non-zero real numbers and suppose that 
\[
\lim_{n\to\infty}\frac{|a_{n+1}|}{|a_n|}=\rho_{\mbox{ratio}}.
\]
Then, if $\rho_{\mbox{ratio}}<1$ the series converges absolutely; if $\rho_{\mbox{ratio}}>1$ or $\rho_{\mbox{ratio}}$ is infinite the series diverges; if $\rho_{\mbox{ratio}}=1$ the test is inconclusive.
\end{theorem}
For the proof see, for example, \cite{HWT}.

\begin{theorem}[The Root Test]
Let $\sum a_n$ be a series of non-zero real numbers and suppose that 
\[
\lim_{n\to\infty}\sqrt[n]{|a_n|}=\rho_{\mbox{root}}.
\]
Then, if $\rho_{\mbox{root}}<1$ the series converges absolutely; if $\rho_{\mbox{root}}>1$ or $\rho_{\mbox{root}}$ is infinite the series diverges; if $\rho_{\mbox{root}}=1$ the test is inconclusive.
\end{theorem}
For the proof see, for example, \cite{HWT}.

Next, we will show the connection between the Ratio Test and the Root Test, constructing the proof that can be used in a Calculus course.
\begin{theorem}[The Relationship Theorem]\label{root-test}
Let $\sum a_n$ be a series of non-zero real numbers and suppose $\rho_{\mbox{ratio}}$ exists (finite or infinite).
Then, $\rho_{\mbox{root}}$ exists and
\[
\rho_{\mbox{root}}=\rho_{\mbox{ratio}}.
\]
The $\rho_{\mbox{root}}$ can be used for convergence verification the same way as $\rho_{\mbox{ratio}}$.
\end{theorem}

\begin{proof}
To simplify notations in this proof, we will occasionally write $\rho$ for $\rho_{\mbox{ratio}}$. First, we will consider finite $\rho$ and then complete the proof for $\rho=\infty$.

Suppose there exists a finite limit: $\lim_{n\to\infty}\frac{|a_{n+1}|}{|a_n|}=\rho \neq 0$. 

By the definition of limit, for every $\epsilon \in (0,\rho)$ there exists $N\in\N$ such that for all $n\ge N$
\[
\rho -\epsilon <\frac{|a_{n+1}|}{|a_n|} < \rho +\epsilon.
\]

In particular,
\[
\rho -\epsilon <\frac{|a_{N+1}|}{|a_N|} < \rho +\epsilon.
\]
Then,
\[
(\rho -\epsilon)\cdot|a_N| < |a_{N+1}| < |a_N|\cdot(\rho +\epsilon)
\]
and 
\begin{equation}\label{ineqal}
(\rho -\epsilon)^k\cdot|a_N| < |a_{N+k}| < |a_N|\cdot(\rho +\epsilon)^k.
\end{equation}
Taking the $(N+k)$-th root, we obtain the estimate:
\[
\sqrt[N+k]{(\rho -\epsilon)^k}\cdot\sqrt[N+k]{|a_N|} < \sqrt[N+k]{|a_{N+k}|} < \sqrt[N+k]{|a_N|}\cdot\sqrt[N+k]{(\rho +\epsilon)^k}.
\]
We can take the limit as $k\to\infty$ in the last inequalities. The following limit is easy to calculate:
\[
\lim_{k\to\infty} \sqrt[N+k]{(\rho \pm \epsilon)^k}\cdot\sqrt[N+k]{|a_N|} = (\rho \pm \epsilon) \cdot 1 = \rho \pm \epsilon
\]
and consequently 
\[
\rho - \epsilon < \lim_{k\to\infty} \sqrt[N+k]{|a_{N+k}|} < \rho+\epsilon
\]
for any positive $\epsilon$.
This implies that 
\[
\rho_{\mbox{root}}= \rho_{\mbox{ratio}}.
\]

\bigskip
If $\rho=0$, we can take arbitrary small $\epsilon \in (0,1)$. Replacing right side of Inequality~(\ref{ineqal}) with
\[
|a_{N+k}| < |a_N|\cdot \epsilon^k
\]
and taking the $(N+k)$-th root, we obtain the estimate:
\[
0\leq\sqrt[N+k]{|a_{N+k}|} < \sqrt[N+k]{|a_N|}\cdot\sqrt[N+k]{\epsilon^k}.
\]
Taking the limit as $k\to\infty$, we can see that
\[
\rho_{\mbox{root}}= 0 = \rho_{\mbox{ratio}}.
\]

\bigskip
If $\rho_{\mbox{ratio}}=\infty$, then for every $M>0$ there exists $N\in\N$ such that for all $n\ge N$
\[
M < \frac{|a_{n+1}|}{|a_n|}
\]
Arguments similar to the above yield the estimate:
\[
M^k \cdot |a_N| < |a_{N+k}|.
\]
Taking the $(N+k)$-th root, we obtain
\[
\sqrt[N+k]{M^k} \cdot \sqrt[N+k]{|a_N|} <\sqrt[N+k]{ |a_{N+k}|}.
\]

Then,
\[ 
\lim_{k\to\infty} \sqrt[N+k]{M^k} \cdot  \sqrt[N+k]{ |a_{N}|} <    \lim_{k\to\infty}\sqrt[N+k]{ |a_{N+k}|},
\]
i.e.,
\[
M <  \lim_{k\to\infty}\sqrt[N+k]{ |a_{N+k}|}
\]
for every $M>0$. This means that $\rho_{\mbox{root}}=\infty$.
\end{proof}

We can think of the tests for convergence in the following ways. The Ratio Test is testing whether the expression 
\begin{equation}\label{max}
\max_{n>L}\left|\frac{a_n}{a_{n-1}}\right|
\end{equation}
is less than $1$ or the expression 
\begin{equation}\label{min}
\min_{n>L}\left|\frac{a_n}{a_{n-1}}\right|
\end{equation}
is greater than $1$ for all $n$ greater than some sufficiently large $L$, while the Root Test is testing whether the geometric mean 
\begin{equation}\label{geom-eqn}
\sqrt[N]{|a_0|\prod_{n=1}^N\left|\frac{a_n}{a_{n-1}}\right|}
\end{equation}  
is less (or greater) than $1$ for all $N$ greater than some sufficiently large $L$. 

We know that the geometric mean~(\ref{geom-eqn}) of a sequence is dominated by the power mean: 
\begin{equation}\label{power-eqn}
\sqrt[p]{\frac{1}{N}\sum_{n=1}^{N}\left|\frac{a_n}{a_{n-1}}\right|^p}, \ \ p\in(0,\infty).
\end{equation}
Note, when $p=1$, the expression~(\ref{power-eqn}) is the arithmetic mean.
\\
The power mean~(\ref{power-eqn}) is dominated by the value of the $\max$~(\ref{max}). This idea leads to the following 
\begin{definition} Let $p\in(0,\infty)$. 
Suppose there exists a finite or infinite limit
\[
 \lim_{N\to\infty}\sqrt[p]{\frac{1}{N}\sum_{n=1}^{N}\left|\frac{a_n}{a_{n-1}}\right|^p}.
\]
Then, we call this limit $\rho_{\mbox{mean}}$.
\end{definition}
The $\rho_{\mbox{mean}}$ can test series convergence with the help of the following  
\begin{theorem}[The Power Mean Tests]\label{power}
Let $\sum_{n=0}^{\infty} a_n$ be a series of non-zero real numbers, let $p\in(0,\infty)$ be the exponent in the definition of $\rho_{\mbox{mean}}$.\\
If $\rho_{\mbox{mean}}<1$, then the series is convergent. 
If $\rho_{\mbox{ratio}}$ exists, then $\rho_{\mbox{mean}}=\rho_{\mbox{root}}=\rho_{\mbox{ratio}}.$ 
\end{theorem}
For the proof of this theorem, we need the following
\begin{lemma}\label{lemma-for-averages}
Let $\left\{b_n\right\}_0^{\infty}\in\R_+$, $I\subset \R$ be a finite or infinite interval, such that $b_0, b_1, ... \in I $ and $f$ be a continuous injective function, $f: I \to \R$. Suppose 
\[
\lim_{N\to\infty} f^{-1}\left(\frac{1}{N}\sum_{n=1}^{N}f(b_n)\right)
\]
exists. Then, for any $L\in\R$
\[
\lim_{N\to\infty} f^{-1}\left(\frac{1}{N}\sum_{n=1}^{N}f(b_n)\right) = \lim_{N\to\infty} f^{-1}\left(\frac{1}{N-L}\sum_{n=L+1}^{N}f(b_n)\right).
\]
\end{lemma}
\begin{proof}[Proof (of the Lemma):]
By continuity of $f^{-1}$
\[
\lim_{N\to\infty}f^{-1}\left(\frac{1}{N-L}\sum_{n=L+1}^{N}f(b_n)\right) =
\lim_{N\to\infty} f^{-1}\left(\frac{1}{N-L}\sum_{n=L+1}^{N}f(b_n) +
\frac{1}{N-L}\sum_{n=1}^{L}f(b_n)\right) .
\]
The latter limit is
\[
\lim_{N\to\infty}f^{-1}\left(\frac{1}{N-L}\sum_{n=1}^{N}f(b_n)\right)=\lim_{N\to\infty}f^{-1}\left(\frac{N}{N-L}\cdot\frac{1}{N}\sum_{n=1}^{N}f(b_n)\right).
\]
 By continuity of $f^{-1}$
\[
\lim_{N\to\infty}f^{-1}\left(\frac{N}{N-L}\cdot\frac{1}{N}\sum_{n=1}^{N}f(b_n)\right)=
\lim_{N\to\infty}f^{-1}\left(\frac{1}{N}\sum_{n=1}^{N}f(b_n)\right).
\]
\end{proof}
Now we can present the
\begin{proof}[Proof (of the Power Mean Tests):]
Suppose $\rho_{\mbox{mean}}<1$. Without loss of generality we can assume that $a_{0}=1$ and $\rho_{\mbox{root}} =  \lim_{N\to\infty} \sqrt[N]{\prod_{n=1}^N\left|\frac{a_n}{a_{n-1}}\right|}$. But
$$
\sqrt[N]{\prod_{n=1}^N\left|\frac{a_n}{a_{n-1}}\right|}
\leq
\sqrt[p]{\frac{1}{N}\sum_{n=1}^{N}\left|\frac{a_n}{a_{n-1}}\right|^p}.
$$
Then, $
 \lim_{N\to\infty}  \sqrt[N]{\prod_{n=1}^N\left|\frac{a_n}{a_{n-1}}\right|}$ is bounded from above by $\rho_{\mbox{mean}}$ that is less than $1$.
Consequently, the series is convergent.

Suppose $\rho_{\mbox{ratio}}$ exists (finite or infinite). Then, for any $\epsilon>0$ there exists sufficiently large $L$ such that  for any $N>L$
$$
\max_{L<n<N}\left|\frac{a_n}{a_{n-1}}\right| < \rho_{\mbox{ratio}}  + \epsilon
$$
and 
$$
 \rho_{\mbox{ratio}}  - \epsilon < \min_{L<n<N}\left|\frac{a_n}{a_{n-1}}\right|.
$$
By Lemma~\ref{lemma-for-averages}, 
\[
\lim_{N\to\infty}\sqrt[p]{\frac{1}{N}\sum_{n=1}^{N}\left|\frac{a_n}{a_{n-1}}\right|^p} = 
\lim_{N\to\infty}\sqrt[p]{\frac{1}{N-L}\sum_{n=L+1}^{N}\left|\frac{a_n}{a_{n-1}}\right|^p}
\]
Also, a power mean is bounded from above and below:
$$
\min_{L<n<N}\left|\frac{a_n}{a_{n-1}}\right|
\leq
\sqrt[p]{\frac{1}{N-L}\sum_{n=L+1}^{N}\left|\frac{a_n}{a_{n-1}}\right|^p}
\leq
\max_{L<n<N}\left|\frac{a_n}{a_{n-1}}\right|.
$$
This implies that $\rho_{\mbox{ratio}}=\rho_{\mbox{mean}}$. Then, the Theorem~\ref{root-test} implies that $\rho_{\mbox{root}}=\rho_{\mbox{mean}}$.
\end{proof}
Similarly to the power mean, we can utilize the generalized $f$-mean (Kolmogorov-Nagumo-de Finetti mean)
\[f^{-1}\left(\frac{1}{N}\sum_{n=1}^{N}f\left|\frac{a_n}{a_{n-1}}\right|\right)\] 
in the following
\begin{definition}
Let $\sum_{n=0}^{\infty} a_n$ be a series of non-zero real numbers, $I\subset \R$ be a finite or infinite interval, such that $\left|\frac{a_1}{a_0}\right|, \left|\frac{a_2}{a_1}\right|, ... \in I $.
Let $f$ be a continuous injective function, $f: I \to \R$. Then $\rho_{\mbox{generalized}}$ is defined in the following way:
$$
\rho_{\mbox{generalized}}=\lim_{N\to\infty} f^{-1}\left(\frac{1}{N}\sum_{n=1}^{N}f\left|\frac{a_n}{a_{n-1}}\right|\right).
$$
\end{definition}
\begin{theorem}[The Generalized $f$-Mean Tests]\label{generalized}
Let $\sum_{n=0}^{\infty} a_n$ be a series of non-zero real numbers, $I\subset \R$ be a finite or infinite interval, such that $\left|\frac{a_1}{a_0}\right|, \left|\frac{a_2}{a_1}\right|, ... \in I $.
Let $f$ be a continuous injective function, $f: I \to \R$,  and suppose that $\rho_{\mbox{generalized}}$ and $\rho_{\mbox{ratio}}$ exist. Then, the value of the $\rho_{\mbox{generalized}}$ is independent of the choice of $f$, in particular, $\rho_{\mbox{generalized}} = \rho_{\mbox{ratio}}$. 

\end{theorem}
\begin{proof}
Similarly to the Theorem~\ref{power}, the conclusion follows from the fact that for any natural numbers $L$ and $N$
$$
\min_{L<n<N}\left|\frac{a_n}{a_{n-1}}\right|
\leq f^{-1}\left(\frac{1}{N-L}\sum_{n=L+1}^{N}f\left|\frac{a_n}{a_{n-1}}\right|\right) \leq
\max_{L<n<N}\left|\frac{a_n}{a_{n-1}}\right|.
$$
\end{proof}

Let us now assume that $p$ is defined on the extended real line and let us denote 
\[
M_p(N):=\sqrt[p]{\frac{1}{N}\sum_{n=1}^{N}\left|\frac{a_n}{a_{n-1}}\right|^p}.
\]
Then, 
\begin{itemize}
\item $\rho_{\mbox{root}} = \lim_{N\to\infty} M_0 (N)$, provided $a_0=1$ (because geometric mean of $N$ real positive numbers is $\lim_{p\to0}M_p(N)$, see \cite{Bu}), 
\item $\rho_{\mbox{mean}} = \lim_{N\to\infty} M_p (N)$ and 
\item $\rho_{\mbox{ratio}}$ is related to the $min (\frac{a_0}{a_1}, \frac{a_1}{a_2},...)= \lim_{N\to\infty} M_{-\infty} (N)$ and to the $max (\frac{a_0}{a_1}, \frac{a_1}{a_2},...)= \lim_{N\to\infty} M_{\infty} (N)$.
\end{itemize}

\begin{remark}\label{instances}
$\rho_{\mbox{ratio}}$, $\rho_{\mbox{root}}$ and $\rho_{\mbox{mean}}$ are the instances of $\rho_{\mbox{generalized}}$.
\end{remark}

\begin{remark}\label{p-comparison}
Note that for $p_1 < p_2$ the Power Mean Test associated with $p_1$ is stronger than the Power Mean Test associated with $p_2$, because 
\[
M_{p_1}(N) < M_{p_2}(N).
\]
\end{remark}

Sometimes, it is easier to find the arithmetic mean of the ratios then to find $\rho_{\mbox{root}}$ or $\rho_{\mbox{ratio}}$. Below we present such 
\begin{example}\label{example} Let us assume $a_0=1$ and define the other terms $a_n$ recursively:
\[
a_{n+1}=\frac{1}{2}\cdot\left( 1-\frac{(-1)^n}{2}  \right)a_n.
\]

The Root Test is computationally difficult. For the Ratio Test, the limit of  $\frac{a_{n+1}}{a_n}$ does not exists.
\\
So, to check the convergence of $\sum_{n=0}^{\infty} a_n$, we will evaluate $\rho_{\mbox{mean}}$.

\[
\rho_{\mbox{mean}} = \lim_{N\to\infty}\frac{1}{N}\sum_{n=1}^{N}\frac{a_n}{a_{n-1}} = \lim_{N\to\infty}\frac{1}{2N}\sum_{n=1}^{N}1-\frac{(-1)^n}{2}  \]
\[
= \lim_{N\to\infty}\frac{1}{2N}\cdot\left( N-\frac{1-(-1)^{N}}{4}  \right) =\frac{1}{2} <1,
\]
and the series $\sum a_n$ is (absolutely) convergent.
\end{example}

\section{``How to be fair''}
The obvious consequence of Theorem~\ref{root-test} is that the radius of convergence of series of functions (and power series in particular) is well defined and can be calculated with either Root or Ratio Tests.

Also, Theorems~\ref{root-test} and (\ref{power}) immediately imply the following statement, which may help us to eliminate unnecessary work, when testing series for convergence. 

\begin{corollary}\label{connections} Given the series $\sum_{n=0}^{\infty}a_n$, the following relations between inconclusiveness of the tests hold:
\begin{itemize}
\item If $\rho_{\mbox{ratio}}=1$, then $\rho_{\mbox{root}}=1$ and $\rho_{\mbox{mean}}=1$.
\item If $\rho_{\mbox{root}}=1$, then either $\rho_{\mbox{ratio}}=1$ or $\rho_{\mbox{ratio}}$ does not exist. 
\item If $\rho_{\mbox{mean}}=1$, then either $\rho_{\mbox{ratio}}=1$ or $\rho_{\mbox{ratio}}$ does not exist.
\item If $\rho_{\mbox{root}}$ does not exist, then $\rho_{\mbox{ratio}}$ does not exist.
\item If $\rho_{\mbox{mean}}$ does not exist, then $\rho_{\mbox{ratio}}$ does not exist.
\end{itemize}
\end{corollary}

"It always does seem to me that I am doing more work than I should do. It is not that I object to the work, mind you; I like work: it fascinates me... But, though I crave for work, I still like to be fair. I do not ask for more than my proper share."  (\cite{JKJ}). 

It would be fair to have Corollary~\ref{connections} in our Calculus course.

\end{document}